%Plain TeX file of the paper: `Curing The Andrews Syndrome'
%by S. B. Ekhad, and D. Zeilberger
%
%begin macros

\baselineskip=14pt
\parskip=10pt
\def\tilde{\widetilde}
\def\Tilde{\char126\relax}

\font\eightrm=cmr8  
\font\eighttt=cmtt8
\magnification=\magstephalf

\parindent=0pt
\overfullrule=0in
%end macros
\bf
\centerline
{
CURING The ANDREWS SYNDROME\footnote{}
{\eightrm Aug. 11, 1996. This paper is accompanied by a Web Page:
{\eighttt http://www.math.temple.edu/\Tilde zeilberg/synd.html}\quad . }
}
\rm
\bigskip
\it
\centerline{
Shalosh B. EKHAD\footnote{$^1$}{\eightrm Supported in parts by the NSF.}
and Doron ZEILBERGER\footnote{$^2$}{\eightrm Supported in part by the NSF.}
           }
\rm
{\bf Abstract:} George Andrews's recent challenge to automated
identity-proving and the WZ method is dealt with. It is argued
that the rivalry between the classical and automated
approaches to hypergeometric sums is beneficial to both.
 
AMS No. 05A15, 33S30, 39A00.
 
 
{\bf PREFACE; In Praise of Bugs}
 
Two summers ago, at the Ann Arbor third meeting on Algebraic Combinatorics
(June 1994), and also at the  Taormina conference in honor of
Adriano Garsia (Aug. 1994), George Andrews described, rather dramatically,
some apparent shortcomings of the {\it WZ method }
(see [PWZ] for a description of
this method). At one moment, he displayed the printout of the
computer-generated proof of a certain hypergeometric
identity, that was outputted by the Maple package {\tt EKHAD} (that now
accompanies [PWZ]). This necessitated
14 transparencies, taped together, which Andrews rolled down on the floor,
producing an obvious comic effect.
 
The identity was ((1.6) of [A1]):
$$
{}_5 F_4 \left (
{{-2i-1,x+2i+2,x-z+1/2,x+i+1,z+i+1} \atop
 {(x+1)/2,x/2+1,2z+2i+2,2x-2z+1}} ; 1 \right ) =0 \quad ,
\eqno(MRR)
$$
where we use the standard hypergeometric notation (see [GKP],[PWZ]).
 
In fact, the printout that was displayed was {\it not even}
the proof of the identity
itself, that was beyond the capabilities of EKHAD at the time, 
but the proof of a certain special case, 
obtained by specializing the parameter $z$.
 
As it turned out, EKHAD's failure to give a complete proof of
$(MRR)$ at the time, as well as the length of the proof of the
special case, were not the algorithm's fault, but were due to
a typo in the Maple coding of the second 
author\footnote{$^3$}{\eightrm 
It escaped notice before because it
occurred in a rarely visited part of the program.}.
Once this typo was corrected, the proof of $(MRR)$ that EKHAD outputs
is rather short (it fits in less than one transparency),
and is produced relatively fast. Both the input and output files can
be viewed in {\tt http://www.math.temple.edu/\Tilde zeilberg/synd.html},
the present paper's Web Page.
 
This was a lucky mistake,
since it lead Andrews to develop an elegant
new technique, that he dubbed `Pfaff's method'([A1-3]). 
Using this new technique, he not only proved $(MRR)$, but in the process
was naturally lead to the {\it discovery and proof} of nineteen additional
{\it brand new} hypergeometric identities. These identities would never
have been discovered had EKHAD been perfect two years ago.
 
After the current version of
EKHAD effortlessly proved identity $(MRR)$, we were sure
that it would just as easily prove all these other identities. This
was implicit in the second author's talk in Herb Wilf's
birthday conference, that took place in Philadelphia, June 12-15, 1996.
When George Andrews asked us whether these new identities
are equally provable by {\tt EKHAD}, in real time, we replied that while we
did not try yet, we were sure that they would pose no problem.
 
Imagine our surprise, and slight embarrassment, when the new, improved,
and debugged, version of EKHAD ran out of memory in each and
every one of Andrews's twenty new identities, except
for two (the original $(4.2)$ (which is $(MRR)$) and
 $(4.20)$ of [A1]). When one specializes
the auxiliary parameters $x$ and $z$ to be specific rational numbers,
then the output is ready after a few seconds. When only one of the
parameters ($x$ or $z$) is specialized, it takes about an hour of CPU
time (on a Sun SPARC), but when both $x$ and $z$ are left alone,
Maple ends with a memory fault (on our system).
 
In this paper we will outline how to get around this problem, and how to
produce absolutely rigorous WZ-proofs for each of these $18$ identities
(and possibly many others yet to be discovered, on which EKHAD might
fail for lack of memory), using commonly available computers. However,
while it is feasible to get such a rigorous proof in each case,
it is very time-consuming, and takes about a week running on our system
on {\tt nice}. On the other hand, the very same method also produces
{\it semirigorous} (see [Z]) 
proofs very fast. Since we know that the identities
are true (even if Andrews did not have a proof), and we know that
we {\it can} find a complete proof, 
if we only were willing to wait a week/identity,
we don't see the point of actually wasting the Temple Math 
Department's computer resources. 
 
Any of our readers who wishes to have a complete
proof is welcome to download our Maple package {\tt SYND}, as well as any
of the input files given in this paper's Web Page, on their own
computer. Alternatively, we would be happy to run the program for
a fee of \$300/identity. Once the fee is paid, 
the program would be run, and once it terminates, the output would be
published in the Web Page, and the name of the donor would be
prominently displayed. This is a good way to honor friends and
relatives, by naming the proof after them. Instructions on how to
order a proof are given in the Web Page mentioned above.
 
In order to demonstrate feasibility, we have run the program
{\tt SYND} on one of Andrews's new identities: $(4.4)$
of [A1]. The input and output files can be gotten from our Web page.
 
{\bf A Short Review of Creative Telescoping and WZ-certification}
 
Let $F(n,k)=F(n,k,c_1, \dots ,c_r)$ be {\it proper hypergeometric} 
(see [PWZ], p. 64 for definition) in all its arguments. The 
Fundamental theorem of [PWZ] guarantees that there exist a
non-negative integer $J$, and
polynomials $a_j(n)=a_j(n,c_1 , \dots , c_r )$, $j=0, \dots , J$ free of
$k$, as well as another proper hypergeometric term
$G(n,k)=G(n,k,c_1, \dots ,c_r)$, such that $G/F$ is a rational function
in all its arguments, and such that:
$$
\sum_{j=0}^{J} a_j(n) F(n+j,k) =G(n,k+1)-G(n,k) \quad .
\eqno(CT)
$$
It follows immediately, by summing w.r.t. $k$, that the definite sum:
$$
A(n)=A(n,c_1, \dots , c_r):=\sum_k F(n,k,c_1, \dots , c_r)
$$
satisfies the linear recurrence equation with polynomial coefficients
$$
\sum_{j=0}^{J} a_j(n) A(n+j) =0 \quad .
\eqno(Recurrence)
$$
Hence, in order to prove a conjectured identity
$$
A(n)=B(n), \quad
$$
where $B(n)$ is a certain explicitly given expression, all we
have to do is verify that $B(n)$ also satisfies the same recurrence, i.e.
$$
\sum_{j=0}^{J} a_j(n) B(n+j) =0 \quad ,
$$
then check that the leading coefficient $a_J(n)$ does not vanish
at positive integers (if it does we just begin at the highest
integer root, and check empirically the finite number of cases before).
Finally we check that $A(j)=B(j)$ for $j=0, \dots , J-1$.
The identity $A(n)=B(n)$ then would follow, by induction, for every
positive integer value of $n$.
 
The recurrence and {\it certificate}, $G(n,k)$, featured in $(CT)$
can be found with the {\it Gosper-Zeilberger} algorithm ([PWZ][GKP]).
 
In this algorithm, in addition to
the unknowns $a_0, \dots , a_J$ one is also looking for
the coefficients, $b_i$, of a certain polynomial:
$$
b(k;n,c_1, \dots, c_r) =\sum_{i=0}^{K} b_i (n, c_1 , \dots , c_r) k^i
$$
such that equation $(6.1.3)$ of [PWZ], p. 108, holds, i.e.
$$
p_2(k)b(k+1)- p_3(k-1) b(k)=p(k) \quad ,
\eqno(G-Z)
$$
for certain explicit polynomials $p(k), p_2(k), p_3(k)$ 
that also depend on $n, c_1, \dots, c_r$, and that
are derivable from the summand $F(n,k)$.
Furthermore $p(k)$ involves the unknown coefficients of the recurrence
$a_j(n)$ {\it linearly}.
 
The next step is to {\it expand} $(G-Z)$ in powers of $k$
and compare the coefficients,
getting a homogeneous system of {\it linear} equations
in the unknowns $a_0, a_1 , \dots ,a_J, b_0, \dots , b_K$.
If the system has no non-trivial solutions,
it means that $J$ was too low, so we try again with $J \leftarrow J+1$. We are
guaranteed to succeed eventually by the Fundamental Theorem (see [PWZ]).
 
Until Andrews's paper [A1], the above was as easily said
as done, and EKHAD had no difficulty in proving any identity of
the above form. Sure, the largest contemporary computer is probably
unable to find the recurrence that
$$
A_n:=\sum_{k=0}^{n} {{n} \choose {k}}^{100000}
$$
satisfies, but for {\it single sum} summations with conjectured
(or already known) explicit expressions, it never failed,
as far as we know. Imagine our chagrin (and George Andrews's justified
glee) at the failure of our program on his new identities.
We should remark that the failure was not {\it theoretical} but
{\it practical}: insufficient memory to handle the large objects
encountered.
 
{\bf The Bottle-Neck that Causes the Andrews Syndrome}
 
As we saw above, the heart of the Gosper-Zeilberger algorithm is the solving a
certain system of linear equations with {\it symbolic} coefficients.
When the summand $F(n,k)$ only depends on $n$ and $k$, and the extra
arguments (parameters) $c_1, c_2, \dots , c_r$ are absent, then
the entries of the coefficient matrix of the system are polynomials
in the single variable $n$. However, when other parameters are present
(in Andrews's case we have $x$ and $z$ in addition to $n$), then
the entries are rather large polynomials of several variables
(in Andrews's case of the three variables $n,x,z$). 
 
It turns out that,
say for $(4.4)$ of [A1], the number of unknowns is $8$:
$a_0, a_1, a_2, a_3$ (the coefficients of the recurrence operator) and
$b_0,b_1,b_2,b_3$ (the coefficients of the polynomial $p(k)$).
In order to find the recurrence operator
$a_0+a_1 N+a_2 N^2 +a_3 N^3$ (where $N$ is the forward shift in
the variable $n$: $Nf(n):=f(n+1)$), and the polynomial
$b(k)=b_0+b_1 k +b_2 k^2+b_3 k^3$, we have to solve a certain linear
system:
$$
{\bf M} \, (a_0, a_1, a_2, a_3, b_0, b_1 , b_2, b_3)^T =0 \quad ,
\eqno(Linear System)
$$
where ${\bf M}$ is an $8 \times 8$ matrix whose entries are certain
(rather large) polynomials of the three variables $n,x,z$. As we said
above, solving this system is beyond Maple's capabilities running
on our system. Our readers are welcome to try EKHAD on any of
Andrews's identities on their system. The input files for EKHAD
(that failed on our system, but may work on bigger ones) are
given in this paper's Web page. We would be happy to hear of
reports of any successful run, and we would announce them on that very
same Web page.

{\bf A general WZ procedure for lowering the order of a recurrence by 1}
 
Suppose $\sum_kF(n,k)=r(n)$ is some identity that we wish to prove. Recall the
WZ paradigm, which says to divide by the right hand side, and to try to
prove instead the identity $f(n):=\sum_k(F(n,k)/r(n))=1$. Well, here is
another possible advantage of this paradigm, which was pointed out by Herb
Wilf: since the recurrence formula that the left side satisfies must have
the solution $f(n)=1$, it follows that the recurrence operator will always
have a right factor of $N-1$, where $N$ is the forward shift in $n$. 
That being
the case we can look for the recurrence that is satisfied by $(N-1)f(n)$, and
it will be {\it of one lower order than the original.} To look for that
recurrence, let's rename the summand $F(n,k)/r(n)$, and call it $F(n,k)$ again.
Then what we want to do is to apply the Zeilberger algorithm to the summand
$F(n+1,k)-F(n,k)$, instead of applying it to $F(n,k)$ itself. We obtain a
recurrence of order 1 less than we would otherwise have found, and we then
need to show that it has only the trivial solution, which we do by
displaying enough initial zero values. This procedure is perfectly general,
and it applies in any situation where the form of the conjectured sum is known.
 
{\bf Getting Around the Andrews Syndrome}
 
Consider again an identity of the form
$$
\sum_k F(n,k,c_1, \dots , c_r)=B(n,c_1, \dots , c_r) \quad .
$$
 
As above, we make the right side $1$, by dividing by $B(n)$.
Renaming $F\leftarrow F/B$, we have to prove an identity
of the form (suppressing the dependence on $c_1, \dots ,c_r$):
$$
A(n):=\sum_k F(n,k)=1 \quad .
$$
It is easy to check whether $A(0)=1$. hence it suffices to prove that
$$
\tilde A(n):=A(n+1)-A(n):=\sum_k [ F(n+1,k)- F(n,k) ]=0 \quad .
\eqno(Efes)
$$
 
In $99\%$ of the cases, one can apply Gosper's algorithm  to the
summand $F(n+1,k)- F(n,k)$
(with respect
to $k$, see [PWZ], chapter 7), and get another hypergeometric term
$G(n,k)$ (which furthermore can be written as $R(n,k)F(n,k)$, where
$R(n,k)$ is  a rational function), such
that
$$
F(n+1,k)- F(n,k) =G(n,k+1)-G(n,k) \quad .
\eqno(WZ)
$$
This is the {\it WZ-equation}. Summing $(WZ)$ with respect to
$k$ immediately yields $(Efes)$. Summing $(WZ)$ with respect to
$n$ yields the so-called {\it companion identity} (see ibid.).
 
However, in the remaining $\%1$ of the cases, that include 
Andrews's twenty identities discussed here, Gosper's algorithm
tells us that no such (hypergeometric) $G(n,k)$ exists. The
more general Gosper-Zeilberger algorithm guarantees, however
(writing $\tilde F (n,k):=F(n+1,k)-F(n,k)$) that there exist
a non-negative integer $J$, and $a_0, \dots , a_J$, 
and $\tilde G(n,k) $ such that
$$
\sum_{j=0}^J a_j(n) \tilde F(n+j,k)=\tilde G(n,k+1) - \tilde G(n,k) \quad.
\eqno(CT')
$$
As before, summing with respect to $k$, we get:
$$
\sum_{j=0}^{J} a_j(n) \tilde A(n+j) =0 \quad .
\eqno(Recurrence')
$$
Once these are found it would
follows that $\tilde A(n) \equiv 0$ provided it is
so for $n=0, \dots , J-1$.
 
Note that the `$J$' of $(Recurrence')$ is one less than that of
$(Recurrence)$, since considering $\tilde A(n)$ rather than
$A(n)$ reduces the order of the recurrence by $1$.
 
While the above process reduces the number of unknowns by $1$, the resulting
symbolic linear system still proves too big for Maple to handle. And now
comes the second {\it twist} on the Gosper-Zeilberger algorithm.
In order to prove that $\tilde A(n) \equiv 0$, we don't actually
have to know what the $a_0, a_1, \dots, a_J$ and the $b_0, \dots , b_K$
are exactly!
All we have to know is that they exist!, plus we have to make sure
that the leading coefficient $a_J(n)$ of the recurrence does not vanish
on any positive integer value of $n$ (or if it does, then we have to know
the highest such value). 
 
The existence of a non-trivial solution of the homogeneous
linear system $(Linear System)$ (when it now applies
to the modified sum as above), is exactly equivalent to the fact
that the determinant of the square matrix ${\bf M}$ vanishes. 
This determinant
is a {\it polynomial} in the variables $n, c_1 , \dots , c_r$.
If we can find {\it a priori} bounds for the degrees in each of
its variables, then plug in enough special cases
(the number of which should be at least equal to
the product of $1$ plus the degrees in each of the variables)
into this determinant, and evaluate these numerous, but fast-to-compute
resulting {\it numerical} determinants, and check whether they are
always $0$, we would have a completely rigorous proof.
 
As hinted in [Z], if you take a non-zero polynomial out of the blue and
plug in random values, it is extremely unlikely that it would be zero.
It such a polynomial yields zero for, say, a hundred different tries,
then it is much more likely that the entire framework of mathematics is flawed
than that this particular polynomial is non-zero. Hence we have a situation
where in order to prove an identity we have the {\it costly}, but feasible,
option of a completely rigorous proof, and an extremely fast and
{\it inexpensive}
way to prove it with probability $1 -\epsilon$. In the package
{\tt SYND}, that comes with this paper, and that can be downloaded
from its Web Page, the procedure {\tt PROOF2} has {\tt Certainty}
as one of its arguments, that can be adjusted by the user. Setting
{\tt Certainty}$:=1$ yields a rigorous, but time-consuming, proof.
Setting it, to say $0.1$ would already yield a $99.9999\%$-sure  proof
very fast.
 
So now we know (either for sure, or almost surely) that the
$a_0(n,c_1, \dots , c_r), \dots ,a_J(n,c_1, \dots , c_r)$, and
the $G(n,k)$ in Eq. $(CT')$ exist, and hence Eq. $(Recurrence')$
holds, for {\it some} $a_0, \dots , a_J$. We don't really care
what they are {\it except}, as mentioned above, we have to make sure
that $a_J(n, c_1 , \dots , c_r)$ does not vanish on positive
integers (and if it does, then to find the largest such).
But this is easy and fast. If $a_J(n , c_1 , \dots , c_r)$ would
have vanished at $n=n_0$, then $(n-n_0)$ would have been a factor,
and it still would have been a factor of
$a_J(n , c_1^{0} , \dots , c_r^{0})$, for any specialization
$c_1=c_1^{0}, \dots , c_r=c_r^{0}$. So all we have to do is
apply procedure {\tt ct} or {\tt zeil} of {\tt EKHAD} to
{\it any} such specialization, and make sure that the polynomial
$a_J$ does not have a factor of the form $(n-n_0)$.
Now the summand only depends
on two symbols: $n$, and $k$, and {\tt ct} runs very fast
(at least on Andrews's identities discussed here).
 
{\bf A Priori Bounds For the Degrees of the Determinant}
 
In order to prove that a polynomial $P(x_1, \dots ,x_m)$ is
the zero polynomial, by plugging in sufficiently many values
for its arguments, we need a priori bounds for its degrees in
each and every one of its arguments. If its degree in $x_i$ is
$d_i$, then checking at all the integer points
$-(d_i+1)/2 \leq x_i \leq (d_i+1)/2 $ ($i=1, \dots , m$),
and verifying that the value of $P$ is zero at each of these
points, would constitute
a rigorous proof that $P$ is the zero polynomial.
 
In the present scenario, the polynomial $P(x_1 , \dots , x_m)$ is
given as a determinant $P:=\det (a_{i,j}(x_1 , \dots , x_m))$,
where $(a_{i,j})$ is a square matrix. To get an upper bound $d_i$(usually
very  dull), for the degree of the determinant in the
variable $x_i$, we replace each entry $a_{i,j}$ by its leading term
with respect to $x_i$. For example $4x^3 z^5 + 3x^2 z^7+5xz$ would be
replaced by $4x^3 x^5$ when the degree in $x$ is sought, and by
$3 x^2 z^7$ when the degree in $z$ is desired. Next, we take the
permanent. Since the permanent of a matrix with monomial entries
can't have any cancellations, the degree in the variable $x_i$ of
the resulting permanent would definitely constitute an upper bound for
$d_i$.

{\bf The Package Synd}
 
The present method (or rather the present twist on the old method) is 
implemented by the Maple package {\tt SYND}, that accompanies
this paper. It can be downloaded from either \break
{\tt http://www.math.temple.edu/\Tilde zeilberg/synd.html}
(download {\tt SYND}), or from \break
{\tt ftp://ftp.math.temple.edu/pub/zeilberg/programs},
or by anonymous {\tt ftp} to {\tt ftp.math.temple.edu}
(login as {\tt anonymous}, password as instructed (your E-mail address),
then {\tt cd pub/zeilberg/programs <CR>}, followed by, 
{\tt get SYND <CR>}.To exit {\tt ftp} you type: {\tt quit <CR>}.
 
Once you have {\tt SYND} on your own computer, make sure that
you have Maple (if you don't, get it!),
and that you are in the directory where {\tt SYND} is. 
Then get into Maple by typing
{\tt maple <CR>}. Once inside Maple, type {\tt read SYND; <CR>}, 
and get the on-line help by typing {\tt ezra();CR}.
 
The two main procedures are {\tt PROOF1} and {\tt PROOF2}.
The former handles sums with one extra parameter, while the
latter handles sums, like those of Andrews[A1], with two extra
parameters ($x$ and $z$). We did not bother to write the analogous
procedures for more auxiliary parameters, but the Maple-literate reader
can do it with no trouble.
 
The function call for {\tt PROOF1} is:
 
{\tt PROOF1(SUMMAND,RHS,k,n,LowerLimit,UpperLimit,a,Certainty);}
 
Here {\tt SUMMAND} is the hypergeometric summand, given either
in terms of factorials, or binomial coefficients, or {\tt rf},
where {\tt rf(a,k)} is the raising factorial
$(a)_k:=a(a+1) \cdots (a+k-1)$; {\tt RHS} is the conjectured
right hand side; {\tt k} is the (single-) summation variable;
{\tt n} is the variable with respect to which the recurrence is
desired; {\tt LowerLimit} and {\tt UpperLimit} are respectively
where the summation starts and ends; {\tt a} is the auxiliary parameter;
Finally, {\tt Certainty} is the rigor-level. Setting
{\tt Certainty:=1;}, would yield a rigorous proof. Setting it
to anything less would give a semi-rigorous proof. Be warned
that {\tt Certainty=0.1} does not mean that the probability that
the identity is true is $\%10$! It is probably true with probability
$\%99.999999$. The meaning of {\tt Certainty} is the fraction of
trials that we attempt, compared to the full intervals in 
{\tt n} and in {\tt a}, that is needed for a rigorous proof. 
 
{\bf Example:} While {\tt EKHAD} can handle the Chu-Vandermonde
identity 
$$
\sum_{k=0}^{n} {{n} \choose {k}} {{a} \choose {k}}= {{n+a} \choose {n}}
\quad ,
$$
in a few seconds, just for the sake of example, here is the function
call that does it in {\tt SYND}. Since {\tt Certainty} is set to $1$,
this would yield a completely rigorous proof.
 
{\tt
PROOF1(binomial(a,k)*binomial(n,k),binomial(a+n,a),k,n,0,n,a,1);
} \quad .

{\tt PROOF2} is exactly as above, except that we have to specify
the two auxiliary parameters: $x$, and $z$.
 
The function call for {\tt PROOF2} is:
 
{\tt PROOF2(SUMMAND,RHS,k,n,LowerLimit,UpperLimit,x,z,Certainty);}
\quad .
 
where {\tt x} and {\tt z} are the auxiliary parameters. For example
the Dixon identity
$$
\sum_{k=-n}^{n} (-1)^k {{a+b} \choose {a+k}} {{a+n} \choose {n+k}}
{{b+n} \choose {b+k}}
= {{(a+b+n)!} \over {a!b!n!}}
\quad ,
$$
is proved, by {\tt SYND}, with the following call:

{\tt
PROOF2((-1)**k*binomial(a+b,a+k)*binomial(a+n,n+k)*binomial(b+n,b+k),
 (a+b+n)!/a!/b!/n!,k,n, -n,n,a,b,1);}
 
(Once again
this is only for the sake of example, since {\tt EKHAD} does it
effortlessly).

{\bf Andrews's 20 Identities}
 
The input files for all the identities of section $4$ of [A1],
as well as the corresponding output files (with {\tt Certainty=0.1}),
can be retrieved from the Web page. As we said above, the only
identity for which we ran the program with {\tt Certainty=1} is
identity $(4.4)$. Readers who have Maple and have idle time to spare,
are welcome to download any or all the input files, change the last
argument of {\tt PROOF2} to $1$ (instead of $0.1$), and run it on their
computer. We request that you notify us, so that we can announce that
the WZ-style rigorous proof of the given identity has been performed,
with due acknowledgement to the computer and the human owner.
\eject
{\bf Lily Yen's Method}
 
Another possibility, that might work,
of getting around the Andrews syndrome is
to adapt Lily Yen's[Y] beautiful approach. Yen found
an a priori bound, $L$ , easily derivable from the identity, such that
if the identity is true for $n=0,1, \dots , L$, then it would be
true for all $n$. Unfortunately, it is so enormous that at present it only
has theoretical interest. However, it seems to us
that an appropriate modification to the situation described
in the present paper, where one has extra parameters, would yield
quite small and practical upper bounds. 
 
The reason for the gargantuan size of Yen's upper bound is
having to bound the largest integer root of the leading
coefficient $a_J(n)$. The upper bound for the order of the
recurrence, $J$, is quite small. Since {\it now} we have
extra parameters, we can rule out the possibility of positive
integer roots as above, by running {\tt ct} of {\tt EKHAD}
on a specialization. It would be interesting and useful to make
this precise, and to implement it.
 
{\bf Peter Paule's Method}
 
Andrews's identities are also interesting from the WZ theory point of
view. The orders of the recurrences outputted by EKHAD ($3$ and $2$) are
two higher than the orders of the minimal recurrences ($1$ and $0$)
satisfied by the right sides.
 
They provide many new non-trivial examples to the phenomenon
described on p. 117 of [PWZ]. This phenomenon, of not getting
the minimal recurrence,
is much more widespread for {\it q-series}, in which Peter Paule[P]
introduced a very useful order-reducing preprocessing device. Paule's
method is also useful for ordinary hypergeometric sums.
 
At present it is not clear how to apply Paule's method to Andrews's
sums, but we suspect that an appropriate generalization will do
the job. Perhaps Andrews's identities are limiting cases of
more general identities that come from multiple sums, on which
there would be some obvious symmetry group with respect to which
one would be able to apply Paule-symmetrization.
 
{\bf Conclusion}
 
Our Wise men, let their memory be blessed said:
`Kinat Sofrim Tarbe Khochma', which, roughly means:
`Rivalry among scholars increases knowledge'. In the present case
the rivalry is between human and machine.
While the meta-mathematical
debate engendered in [Z] and [A4] is unlikely to be resolved in our
time, mathematics proper does benefit. 
 
The second moral is that {\it mistakes} are crucial for progress.
According to the current dogma in biology, we humans would have
still been amoebas if not for a series of lucky mistakes in
biological transmission of information. The same can be said 
for science, and even for mathematics. If the previous version
of {\tt EKHAD} did not contain a bug, George Andrews would have
obtained the proof of $(MRR)$ right away, and would never have needed
to find another proof, that lead to twenty new 
beautiful identities.
While the WZ method is also capable of discovering
(and proving at the same time) new identities out of old ones,
we do not know at present how to rediscover, naturally, Andrews's
identities from scratch.
 
Bruno Salvy has told us that the question of deciding how many
random checks of the vanishing of a polynomial are needed in order
to deduce its nullity, within a prescribed margin of error, has been
dealt with. See his message[B] that can be accessed from this
paper's Web page.
 
One last point. 
The present method is easy to parallelize.
This is because it boils down to checking many special instances,
that can be done independently of each other. If the need would arise 
in the future to
prove a huge identity (for example because it would
imply the Riemann Hypothesis), then it would be feasible to have
an International collaboration of many computers.
 
{\bf Acknowledgement}
 
We wish to thank George Andrews and Herb Wilf for many stimulating
discussions.
 
{\bf REFERENCES}
 
[A1] George Andrews,
{\it Pfaff's Method (I): The Mills-Robbins-Rumsey determinant},
Discrete Mathematics (Garsia Festschrifft issue), to appear.
 
[A2] George Andrews,
{\it Pfaff's Method (II): diverse applications}, to appear.
 
[A3] George Andrews,
{\it Pfaff's Method (III): comparison with the W-Z method},
Elect. J. of Combinatorics {\bf 3(2)}(1996) [Foata Festschrifft], R21.
 
[A4] George Andrews,
 {\it The death of proof? Semi-rigorous mathematics? you've got to
be kidding!}, Math. Intell. {\bf 16}, no. 4, 16-18 (Fall 1994).
 
[B] Bruno Salvy, {\it Message to Doron Zeilberger}, published,
by permission, in this paper's Web page
{\tt http://www.math.temple.edu/\Tilde zeilberg/synd.html}.
 
[GKP] Ronald Graham, Donald Knuth, and Oren Patashnik,
{\it ``Concrete Mathematics''}, second edition, Addison-Wesley,
Reading, 1994.
 
[P] Peter Paule,
{\it  Short and easy computer proofs of the Rogers-Ramanujan identities
and of identities of similar type},
Elect. J. of Combinatorics {\bf 1}(1994) R10.
 
[PWZ] Marko Petkovsek, Herbert Wilf, and Doron Zeilberger,
{\it ``A=B''}, AK Peters, Wellesley, (1996).
{\tt Web Page: http://central.cis.upenn.edu/\Tilde wilf/AeqB.html}\quad .
 
[Y] Lily Yen, {\it 
A two-line algorithm for proving terminating hypergeometric identities},
J. Math. Anal. Appl. {\bf 198}(1996), 856-878.
 
[Z] Doron Zeilberger,
 {\it Theorems for a price: Tomorrow's semi-rigorous mathematical culture},
Notices of the Amer. Math. Soc. {\bf 40 \# 8}, 978-981 (Oct. 1993).
Reprinted: Math. Intell. {\bf 16}, no. 4, 11-14 (Fall 1994).
 
\medskip
-------------------
\medskip
 
Shalosh B. Ekhad, Department of Mathematics, Temple University,
Philadelphia, PA 19122. \break
{\tt E-mail: ekhad@math.temple.edu} \quad ; \quad
{\tt WWW: http://www.math.temple.edu/\Tilde ekhad} \quad .
 
Doron Zeilberger, Department of Mathematics, Temple University,
Philadelphia, PA 19122. \break
{\tt E-mail: zeilberg@math.temple.edu} \quad ; \quad
{\tt WWW: http://www.math.temple.edu/\Tilde zeilberg} \quad .
 
\bye